\documentclass[12pt]{article}
\usepackage{latexsym}
\usepackage{times}
\usepackage{graphicx,color}
\usepackage{amsmath,amssymb, amsfonts}
\catcode`@=11
\renewcommand{\@begintheorem}[2]{\it \trivlist            
      \item[\hskip \labelsep{\bf #1\ #2{\rm :}}]}         
\renewcommand{\@opargbegintheorem}[3]{\it \trivlist       
      \item[\hskip \labelsep{\bf #1\ #2\ {\rm (#3)\/:}}]}
\def\@sect#1#2#3#4#5#6[#7]#8{\ifnum #2>\c@secnumdepth
     \def\@svsec{}\else 
     \refstepcounter{#1}\edef\@svsec{\csname the#1\endcsname{.}\hskip 1em }\fi
     \@tempskipa #5\relax
      \ifdim \@tempskipa>\z@ 
        \begingroup #6\relax
          \@hangfrom{\hskip #3\relax\@svsec}{\interlinepenalty \@M #8\par}
        \endgroup
       \csname #1mark\endcsname{#7}\addcontentsline
         {toc}{#1}{\ifnum #2>\c@secnumdepth \else
                      \protect\numberline{\csname the#1\endcsname}\fi
                    #7}\else
        \def\@svsechd{#6\hskip #3\@svsec #8\csname #1mark\endcsname
                      {#7}\addcontentsline
                           {toc}{#1}{\ifnum #2>\c@secnumdepth \else
                             \protect\numberline{\csname the#1\endcsname}\fi
                       #7}}\fi
     \@xsect{#5}}

\@addtoreset{equation}{section}
\newcommand{\Delete}[1]{}
\newcommand{\pend}{\hspace*{\fill} $\Box$}

\newtheorem{lemma}{Lemma}[section]
\newtheorem{theorem}[lemma]{Theorem}

\newtheorem{remark}{Remark}

\begin{document}

\title{Compression of M${}^\natural$-convex Functions --- Flag Matroids and Valuated Permutohedra}

\author{Satoru FUJISHIGE\footnote{
Research Institute for Mathematical Sciences, 
Kyoto University, Kyoto 606-8502, Japan. 
E-mail: fujishig@kurims.kyoto-u.ac.jp}\ \ 
 and\ \ Hiroshi HIRAI\footnote{
Department of Mathematical Informatics,
Graduate School of Information Science and Technology,
The University of Tokyo, Tokyo, 113-8656, Japan.
E-mail: hirai@mist.i.u-tokyo.ac.jp
}}

\date{May 26, 2020}

\maketitle

\begin{abstract}

Murota (1998) and Murota and Shioura (1999) introduced concepts of 
M-convex function and M${}^\natural$-convex function as discrete convex 
functions, which are generalizations of valuated matroids due to Dress and 
Wenzel (1992).
In the present paper we consider a new operation defined by a convolution of sections 
of an M${}^\natural$-convex function that transforms the given  
M${}^\natural$-convex function to an M-convex function, which we call a 
{\it compression} of an M${}^\natural$-convex function.
For the class of valuated generalized matroids, which are 
special M${}^\natural$-convex functions, the compression induces 
a {\it valuated permutohedron} together with 
a decomposition of the valuated generalized matroid into 
{\it flag-matroid strips}, 
each corresponding to a maximal linearity domain of the induced valuated 
permutohedron. 
We examine the details of the structure of flag-matroid strips and the induced 
valuated permutohedron by means of discrete convex analysis of Murota.

\end{abstract}


\section{Introduction}
\label{sec:1}

Murota (1998) and Murota and Shio\-ura (1999) 
introduced the concepts of M-convex function \cite{Murota98b} and 
M${}^\natural$-convex function \cite{Murota+Shioura99}, as discrete convex 
functions. Their original ideas can be traced back 
to Dress and Wenzel's valuated matroids \cite{Dress+Wenzel92} introduced 
in 1992. 
See \cite{Murota03a,Murota2016,Murota2019} for details about 
the theory of discrete convex analysis
and its applications developed by Murota and others 
(also see \cite[Chapter VII]{Fuji05}).

In the present paper we consider a new operation  defined by 
a convolution of sections 
of an M${}^\natural$-convex function that transforms the given  
M${}^\natural$-convex function to an M-convex function, which we call a 
{\it compression} of an M${}^\natural$-convex function.
For the class of valuated generalized matroids, which is 
a special class of M${}^\natural$-convex functions, the compression induces 
a valuated permutohedron together with 
a decomposition of the valuated generalized matroid into flag-matroid strips, 
each corresponding to a maximal linearity domain of the valuated permutohedron. 
We examine the details of the structure of flag-matroid strips and the induced 
valuated permutohedron.
We investigate the structures of the strip decomposition of
valuated generalized-matroids, special 
M${}^\natural$-convex functions by means of discrete convex analysis 
of Murota (\cite{Murota03a}).
The strip decomposition of a valuated generalized-matroid
uniquely determines a valuated permutohedron, identified with 
a special M-convex function. 

The present paper is organized as follows. In Section~\ref{sec:2} we give 
some definitions and preliminaries about (i) submodular/supermodular functions
and related polyhedra such as base polyhedra, submodular/supermodular 
polyhedra, and generalized polymatroids, and (ii) M-/M${}^\natural$-convex 
functions and L-/L${}^\natural$-convex functions. 
In Section~\ref{sec:compression} we introduce a new operation called 
a {\it compression} of M${}^\natural$-convex functions, which leads us to  
the concepts of flag-matroid strips and a strip decomposition of 
M${}^\natural$-convex functions in Section~\ref{sec:strip}. 
In Section~\ref{sec:VGM} we consider valuated generalized matroids, which are 
special M${}^\natural$-convex functions, and examine implications of 
our results in valuated generalized matroids. 
The strip decomposition of a valuated generalized matroid gives 
a collection of strips of
flag matroids~\cite{Borovik+Gelfand+White03}, each inducing 
a sub-permutohedra.   
The compression of a valuated generalized matroid induces a valuated 
permutohedron
whose maximal linearity domains corresponds to flag-matroid strips 
of the valuated generalized matroid.
Section~\ref{sec:conclusion} gives some concluding remarks.

\section{Definitions and Preliminaries}\label{sec:2}

Let $E=[n](=\{1,\cdots,n\})$ for a positive integer $n>1$. 
For any $x\in\mathbb{R}^E$ and $X\subseteq E$ define $x(X)=\sum_{e\in X}x(e)$, 
where $x(\emptyset)=0$. 
For any subset $X\subseteq E$ its {\it characteristic vector} $\chi_X$ 
in $\mathbb{R}^E$ is defined by $\chi_X(e)=1$ if $e\in X$ and $\chi_X(e)=0$
if $e\in E\setminus X$. We also write $\chi_e$ 
instead of $\chi_{\{e\}}$ for $e\in E$.

\subsection{Basics of submodular/supermodular functions} 

A function $f: 2^E\to\mathbb{R}$ is called a {\it submodular function}
if it satisfies
\begin{equation}\label{eq:1}
 f(X)+f(Y)\ge f(X\cup Y)+f(X\cap Y)\qquad (\forall X, Y\subseteq E).
\end{equation}
We assume that $f(\emptyset)=0$ for any set 
function $f: 2^E\to\mathbb{R}$  in the
sequel. A negative of a submodular function is called a {\it supermodular 
function}. (Also see \cite{Edmonds70,Fuji05}.)

For a submodular function $f: 2^E\to\mathbb{R}$ define 
\begin{equation}\label{eq:2}
 {\rm P}(f)=\{x\in\mathbb{R}^E\mid \forall X\subseteq E: x(X)\le f(X)\},
\end{equation}
which is called the {\it submodular polyhedron} associated with submodular
function $f$. Also define
\begin{equation}\label{eq:3}
 {\rm B}(f)=\{x\in{\rm P}(f)\mid x(E)=f(E)\},
\end{equation}
which is called the {\it base polyhedron} associated with submodular
function $f$. 
As is well known (see \cite{Fuji05}), the base polyhedron ${\rm B}(f)$ 
is always nonempty (and is a face of ${\rm P}(f)$). 

For a supermodular function $g: 2^E\to\mathbb{R}$ we define in a dual manner
the {\it supermodular polyhedron}
\begin{equation}\label{eq:4}
 {\rm P}(g)=\{x\in\mathbb{R}^E\mid \forall X\subseteq E: x(X)\ge g(X)\}
\end{equation}
and the {\it base polyhedron}
\begin{equation}\label{eq:5}
 {\rm B}(g)=\{x\in{\rm P}(g)\mid x(E)=g(E)\}.
\end{equation}
For a submodular function $f: 2^E\to\mathbb{R}$ define a supermodular 
function $f^{\#}: 2^E\to\mathbb{R}$ by
\begin{equation}\label{eq:6}
 f^{\#}(X)=f(E)-f(E\setminus X)\qquad (\forall X\subseteq E),
\end{equation}
which is called the {\it dual supermodular function} of $f$.
Then we have ${\rm B}(f)={\rm B}(f^{\#})$. 
For more details about submodular/supermodular functions and associated 
polyhedra see \cite[Chapter II]{Fuji05}, 
where  submodular/supermodular functions 
defined on distributive lattices ${\mathcal D}\subseteq 2^E$ are also 
investigated and their base polyhedra are unbounded 
unless $\mathcal{D}=2^E$.

Define $Q_{\mathbb{Z}}=Q\cap \mathbb{Z}^E$ for any set  
$Q\subseteq \mathbb{R}^E$. 
Also denote by ${\rm Conv}(Q)$ the convex hull of $Q$ in $\mathbb{R}^E$.
When ${\rm Conv}(Q_{\mathbb{Z}})=Q$, we identify $Q$ 
with $Q_{\mathbb{Z}}$. 

When a submodular function $f: 2^E\to\mathbb{R}$ is integer-valued, 
its submodular polyhedron and base polyhedron are integral (i.e., every 
vertex of the polyhedra is an integral vector). 
Moreover, we have
\begin{equation}\label{eq:7}
 {\rm Conv}({\rm P}(f)_{\mathbb{Z}})={\rm P}(f), \qquad
 {\rm Conv}({\rm B}(f)_{\mathbb{Z}})={\rm B}(f).
\end{equation}
Most of the following arguments are valid even if we regard $\mathbb{Z}$ 
in place of $\mathbb{R}$ as the underlying totally ordered additive group. 
When $f$ is integer-valued, 
we call ${\rm P}(f)_{\mathbb{Z}}$ and ${\rm B}(f)_{\mathbb{Z}}$ 
the submodular polyhedron and the base polyhedron, respectively,
 associated with $f$ as well.
(This is the approach taken in \cite{Fuji05}, indeed.)

For a submodular function $f: 2^E\to\mathbb{R}$ and a nonempty proper 
subset $A$ of $E$
define functions $f^A: 2^A\to\mathbb{R}$ and 
$f_A: 2^{E\setminus A}\to\mathbb{R}$ by
\begin{equation}\label{eq:9}
 f^A(X)=f(X)\quad (\forall X\subseteq A),\quad 
 f_A(X)=f(X\cup A)-f(A)\quad (\forall X\subseteq E\setminus A).
\end{equation}
We also define $f^E=f_\emptyset = f$. 
We call $f^A$ the {\it restriction} of $f$ on $A$ and $f_A$ the 
{\it contraction} of $f$ by $A$.
Similarly we define the restriction and contraction for supermodular functions.

\subsection{Permutohedra and sub-permutohedra}\label{sec:sub-p}

For a permutation $\pi: [n]\to[n]$ we identify $\pi$ with 
the {\it permutation vector} 
$(\pi(1),\cdots,\pi(n))$ in $\mathbb{Z}^n$, which we denote 
by $v_\pi$ (or simply by $\pi$ 
when there is no possibility of confusion). 
Suppose that an integer-valued submodular function $f: 2^E\to\mathbb{Z}$
is given by $f(X)=\sum_{i=1}^{|X|}(n-i+1)$ for all $X\subseteq E=[n]$. 
Then the base polyhedron 
${\rm B}(f)$ has the $n!$ extreme points, each being a permutation vector 
identified with 
a permutation of $[n]$, which is called the {\it permutohedron} 
(or {\it permutahedron}) in $\mathbb{R}^E$.

For any permutation $\pi$ of $[n]$ we have a unique complete flag
\begin{equation}\label{eq:sf1}
  \mathcal{F}: F_1\subset \cdots \subset F_n=[n]
\end{equation}
such that for each $i\in[n]$ $F_i$ is the set of the first $i$ elements of
$(\pi(1),\cdots,\pi(n))$, i.e.,
\begin{enumerate}
\item $|F_i|=i$ for each $i=1,\cdots,n$ and 
\item $\displaystyle{\sum_{i=1}^n\chi_{F_i}=v_\pi}$, where $\chi_{F_i}$ 
is the characteristic vector of a set $F_i\subseteq [n]$.
\end{enumerate}
Denote the flag $\mathcal{F}$ in (\ref{eq:sf1}) by 
$\mathcal{F}^\pi: F^\pi_1\subset \cdots\subset F^\pi_n$.

Let us consider a polyhedron $P$ satisfying the following two: 
\begin{itemize}
\item[(P1)] $P$ is the convex hull of a set of some permutation vectors.
\item[(P2)] $P$ is a base polyhedron. 
\end{itemize}
We call such a polyhedron $P$ a {\it sub-permutohedron}.\footnote{
We may call the sub-permutohedron a {\sl permutohedron}, and 
an ordinary permutohedron a {\sl complete permutohedron}, but we resist
the temptation.}
A sub-permutohedron is precisely the Coxeter matroid polytope
associated with a flag matroid of complete flag (see 
\cite{Borovik+Gelfand+White03} and the discussions to be made 
in Section~\ref{sec:VGM}).
A  recent interesting appearance of a sub-permutohedron
is from the theory of  Bruhat order (\cite{BjoernerBrenti2005}), 
due to Tsukerman and Williams \cite{TsukermanWilliams2015}, that 
every Bruhat interval polytope is a sub-permutohedron.

\subsection{Base polyhedra, generalized polymatroids, and strong maps}

Suppose that a submodular function $f: 2^E\to\mathbb{R}$ and 
a supermodular function $g: 2^E\to\mathbb{R}$ with 
$f(\emptyset)=g(\emptyset)=0$ satisfy the following inequalities
\begin{equation}\label{eq:10}
  f(X)-g(Y)\ge f(X\setminus Y)-g(Y\setminus X)\quad (\forall X, Y\subseteq E).
\end{equation}
Then the polyhedron
\begin{equation}\label{eq:11}
 {\rm P}(f,g)
 =\{x\in\mathbb{R}^E\mid \forall X\subseteq E: g(X)\le x(X)\le f(X)\}
\end{equation}
is called a {\it generalized polymatroid} (\cite{Frank81b,Hassin82}).
There exists a one-to-one correspondence between base polyhedra and 
generalized polymatroids up to translation along a coordinate axis as follows
(see \cite[Figure 3.7]{Fuji05}).

\begin{theorem}[\cite{Fuji84a,Fuji05}]\label{th:b-gp}
For the base polyhedron 
${\rm B}(f)$ associated with a submodular function $f: 2^E\to\mathbb{R}$
 the projection of ${\rm B}(f)$ along an axis $e \in  E$ on the 
coordinate subspace given by the 
hyperplane $x(e) = 0$ is a generalized polymatroid 
${\rm P}(f^\prime,g^\prime)$ in $\mathbb{R}^{E^\prime}$ with 
$E^\prime = E - \{e\}$, where 
$f^\prime$ is the restriction of $f$ on ${E'}$ and 
$g^\prime$ is the restriction of $f^\#$ on ${E'}$.

Conversely, every generalized polymatroid in $\mathbb{R}^{E^\prime}$ 
is obtained in this way. 
\end{theorem}

When a generalized polymatroid is a convex hull of $\{0,1\}$-valued 
points (vertices), then 
it is called a {\it generalized-matroid polytope},
which can be identified with the family ${\mathcal G}$ of subsets 
$X\subseteq E$ such that $\chi_X$ are vertices of 
the generalized-matroid polytope. 
The family $\mathcal G$
is called a {\it generalized matroid} (see \cite{Frank+Tardos88}).
 
An ordered pair $(f_1,f_2)$ of submodular functions 
$f_i: 2^E\to\mathbb{R}$ $(i=1,2)$ is called a {\it weak map} if we have
\begin{equation}\label{eq:12}
 {\rm P}(f_2)\subseteq {\rm P}(f_1).
\end{equation}
Moreover, the ordered pair $(f_1,f_2)$ is called a {\it strong map} 
if we have
\begin{equation}\label{eq:13}
 {\rm P}((f_2)_X)\subseteq {\rm P}((f_1)_X)
 \qquad (\forall X\subset E),
\end{equation}
i.e., every ordered  pair $((f_1)_X,(f_2)_X)$ of contractions of $f_1$ and 
$f_2$ by $X\subset E$ is a weak map. 
The concept of strong map was originally considered for matroids 
(see \cite{Oxley92,Schrijver03,Welsh76}), 
and we adopt it to any submodular functions (or submodular systems).

We also have the following theorem.

\begin{theorem}\label{th:sm-gp}
${\rm P}(f,g)$ is a generalized polymatroid if and only if 
$(f,g^{\#})$ is a strong map.
\medskip\\
{\rm (Proof)
Relation (\ref{eq:13}) is equivalent to the following inequalities
\begin{equation}\label{eq:13a}
 f_2(Z\cup X)-f_2(X)\le f_1(Z\cup X)-f_1(X)
 \qquad (X\subset E,\ Z\subseteq E\setminus X).
\end{equation}
Putting $W=E\setminus(Z\cup X)$, (\ref{eq:13a}) is rewritten 
in terms of the dual supermodular function $f_2^\#$ of $f_2$ as
\begin{equation}\label{eq:13b}
 f_2^\#(Z\cup W)-f_2^\#(W)\le f_1(Z\cup X)-f_1(X)
 \qquad (X\subset E,\ Z\subseteq E\setminus X)
\end{equation}
with $W=E\setminus(Z\cup X)$. Because of the supermodularity of $f_2^\#$ 
(\ref{eq:13b}) is equivalent to
\begin{equation}\label{eq:13c}
 f_2^\#(Z\cup W)-f_2^\#(W)\le f_1(Z\cup X)-f_1(X)
\end{equation}
for all $X, W, Z\subset E$ with $X\cap W=X\cap Z=Z\cap W=\emptyset$, 
which is equivalent to (\ref{eq:10}). 
\pend}
\end{theorem}

A sequence of submodular functions $f_1,\cdots, f_p: 2^E\to\mathbb{R}$ 
is called a {\it strong map sequence} if for each $i=1,\cdots,$\ $p$$-$$1$ 
the pair $(f_{i+1},f_i)$ is a strong map. It follows from 
Theorem~\ref{th:sm-gp} that 
\begin{itemize}
\item[(F1)] Given a strong map sequence 
$f_1,\cdots, f_p: 2^E\to\mathbb{R}$, we have a sequence of generalized
polymatroids ${\rm P}(f_{i+1},(f_i)^{\#})$ for $i=1,\cdots,p-1$.
\end{itemize}

We also have the following.
\begin{itemize}
\item[(F2)] 
For a generalized polymatroid ${\rm P}(f,g)$ and $\alpha\in\mathbb{R}$
such that $g(E)\le \alpha \le f(E)$ the intersection of  ${\rm P}(f,g)$ and
the hyperplane $x(E)=\alpha$ is a base polyhedron (we call such a base polyhedron 
a {\it section} of ${\rm P}(f,g)$ and denote it by ${\rm P}(f,g)_{(\alpha)}$). 
When ${\rm P}(f,g)$ is integral and $\alpha$ is an integer, the section 
${\rm P}(f,g)_{(\alpha)}$ is integral. 
Moreover, letting ${\rm B}(f')$ be a section of ${\rm P}(f,g)$ 
with a submodular function $f'$, 
we have a strong map sequence $g^{\#}, f', f$, i.e., ${\rm P}(f',g)$
and ${\rm P}(f,(f')^{\#})$ are generalized polymatroids.
\end{itemize}

A strong map sequence $f_1,\cdots, f_p: 2^E\to\mathbb{Z}$ with 
$f_i$ $(i=1,\cdots, p)$ being matroid rank functions is called a 
{\it flag matroid} (see \cite{Borovik+Gelfand+White03}).

\subsection{M${}^\natural$-convex functions and L${}^\natural$-convex 
functions}\label{sec:M-natural}

\begin{figure}[h]
 \begin{center}
  \leavevmode
  \includegraphics[width=10cm,height=7cm,clip]{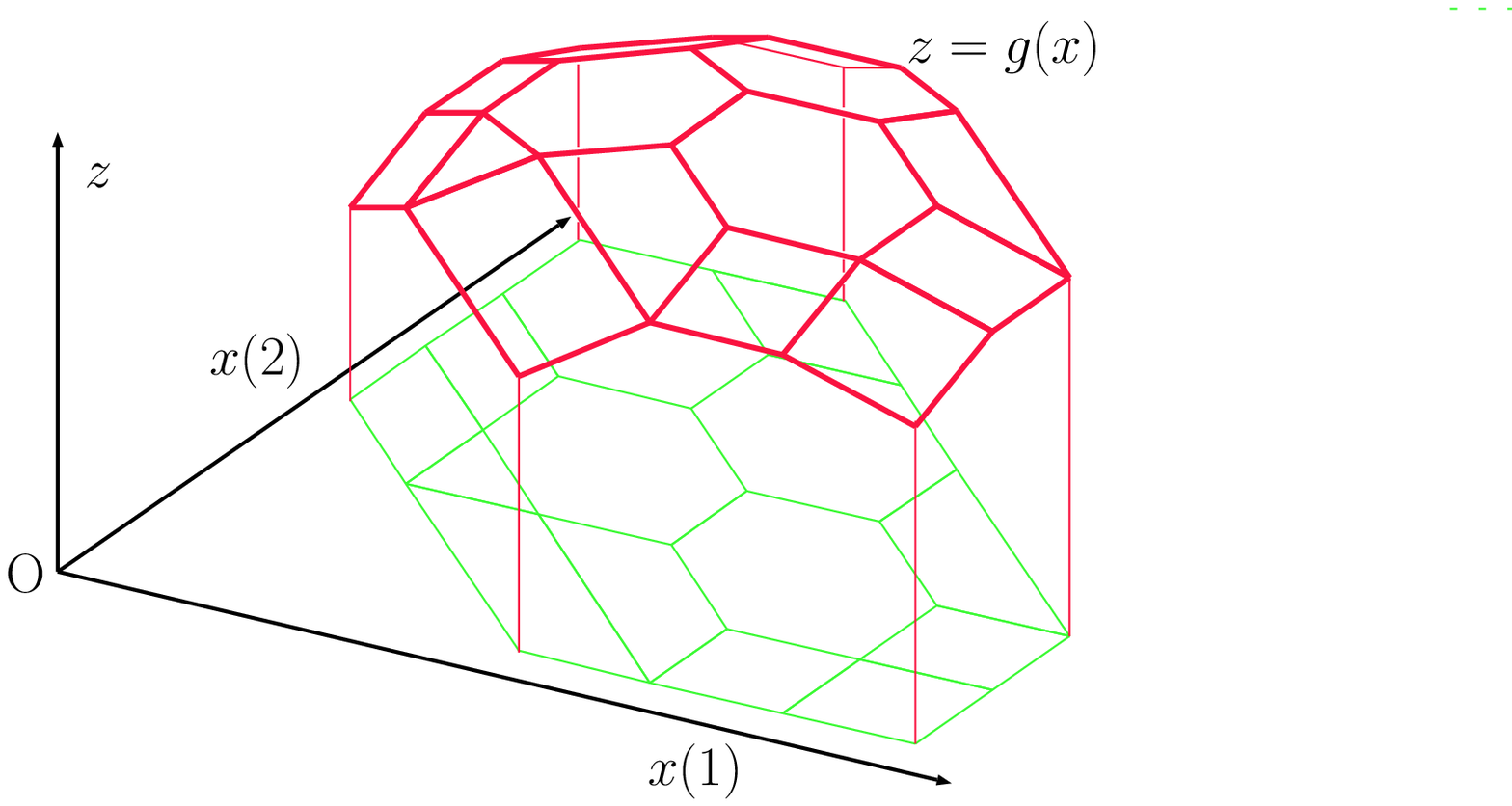}
 \end{center}
\caption{An M${}^\natural$-concave function $g$ (\cite[Fig.~17.4]{Fuji05}).}
\label{fig:M-natural}
\end{figure}

Let $f: \mathbb{R}^E\to\mathbb{R}\cup\{+\infty\}$ be a polyhedral 
convex function such that 
\begin{enumerate}
\item its effective domain, ${\rm dom}f\equiv\{x\in\mathbb{R}^E\mid 
f(x)<+\infty\}$, is a generalized polymatroid (hence nonempty) and
\item every linearity domain of $f$ is a generalized polymatroid,
\end{enumerate}
where a {\it linearity domain} (or {\it affinity domain}) 
of $f$ is ${\rm Arg}\min(f-h)$ for a  
linear function $h(x)=\langle z,x\rangle(\equiv\sum_{e\in E}z(e)x(e))$ 
with some $z\in(\mathbb{R}^E)^*$.
Then $f$ is called an {\it M${}^\natural$-convex 
function}\footnote{Here we employ another equivalent definition of 
M${}^\natural$-convexity, instead of the original definition by means of
the exchange axiom introduced by Murota and Shioura 
\cite{Murota+Shioura99,Murota03a}.}, which is due to
Murota and Shioura \cite{Murota+Shioura99,Murota03a} 
(also see \cite[Section~17]{Fuji05}). 
The negative of an M${}^\natural$-convex function is called an {\it 
M${}^\natural$-concave function} (see Figure~\ref{fig:M-natural}). 

When ${\rm dom}f$ of an 
M${}^\natural$-convex function $f$ is a base polyhedron, $f$ is called an
{\it M-convex function} (see \cite{Murota03a}).\footnote{
M-convex functions were introduced earlier than M${}^\natural$-convex 
functions by Murota (see \cite{Murota98b,Murota03a}). There is a one-to-one 
correspondence between M-convex functions and M${}^\natural$-convex convex
functions because of Theorem~\ref{th:b-gp}.}
Any concept related to M-/M${}^\natural$-concave functions is defined 
in a natural way from that defined for M-/M${}^\natural$-convex functions.

Now let $f: \mathbb{R}^E\to\mathbb{R}\cup\{+\infty\}$ be 
an M${}^\natural$-convex function satisfying  
the following (a) and (b):
\begin{itemize}
\item[(a)] The effective domain of $f$
is an {\sl integral} generalized polymatroid. 
\item[(b)] Every linearity domain of $f$ is also
an {\sl integral} generalized polymatroid. 
\end{itemize}
Such an M${}^\natural$-convex function 
$f$ can be identified with $f$ being restricted on the integer lattice 
$\mathbb{Z}^E$.  
So we can consider an M${}^\natural$-convex function
$f: \mathbb{Z}^E\to\mathbb{R}\cup\{+\infty\}$. 
An integer-valued M-concave function 
$g:\mathbb{Z}^E\to\mathbb{Z}\cup\{+\infty\}$ 
with its effective domain being a matroid base polytope coincides  
with a {\it valuated matroid} due to Dress and Wenzel~\cite{Dress+Wenzel92}.
Also, if an M${}^\natural$-convex function 
$f: \mathbb{Z}^E\to\mathbb{R}\cup\{+\infty\}$ has an effective domain 
${\rm dom}(f)$ whose convex hull ${\rm Conv}({\rm dom}(f))$ is 
 a permutohedron, we call $f$ a {\it valuated permutohedron}.

For any M${}^\natural$-convex function 
$f: \mathbb{R}^E\to\mathbb{R}\cup\{+\infty\}$ define the 
{\it Legendre-Fenchel transform}
(or  {\it convex conjugate}) of $f$ by
\begin{equation}\label{eq:L1}
 f^\bullet(y)=\sup\{\langle y,x\rangle - f(x)\mid x\in\mathbb{R}^E\} 
 \qquad (y\in(\mathbb{R}^E)^*),
\end{equation}
where $\langle y,x\rangle=\sum_{e\in E}y(e)x(e)$. 
The function $f^\bullet$ is called an {\it L${}^\natural$-convex function} 
(\cite{Murota03a}),
which is equivalent to submodular integrally convex function due to Favati 
and Tardella~\cite{FavatiTardella1990}. 
The original $f$ is recovered from $f^\bullet$ by taking another 
Legendre-Fenchel transform as follows.
\begin{equation}\label{eq:L2}
 f(x)=\sup\{\langle y,x\rangle - f^\bullet(y)\mid y\in(\mathbb{R}^E)^*\} 
 \qquad (x\in\mathbb{R}^E).
\end{equation}
(See \cite{Murota03a,Murota2016}.) 
Hence there exists a one-to-one correspondence between 
M${}^\natural$-convex functions and L${}^\natural$-convex functions. 
Furthermore, Murota \cite{Murota03a} showed the integrality property that 
(\ref{eq:L1}) and (\ref{eq:L2}) with $\mathbb{R}$ being replaced 
by $\mathbb{Z}$ 
hold for any M${}^\natural$-convex function 
$f: \mathbb{Z}^E\to\mathbb{Z}\cup\{+\infty\}$.
When $f$ is an M-convex function, its Legendre-Fenchel transform is 
what is called an {\it L-convex function} (\cite{Murota03a}).

For any $x\in\mathbb{R}^E$ the {\it subdifferential} of $f$ at $x$,
denoted by $\partial f(x)$, is defined by
\begin{equation}\label{eq:L3}
\partial f(x)=\{w\in(\mathbb{R}^E)^*
    \mid \forall z\in\mathbb{R}^E: f(z)\ge f(x)+\langle w,z-x\rangle\}.
\end{equation}
The subdifferential of $f^\bullet$ is defined similarly as
\begin{equation}\label{eq:L3a}
\partial f^\bullet(w)=\{x\in\mathbb{R}^E
    \mid \forall y\in(\mathbb{R}^E)^*: 
           f^\bullet(y)\ge f^\bullet(w)+\langle y-w,x\rangle\}.
\end{equation}
Then we have the following.

\begin{lemma}\label{lem:sbd1}
We have $w\in \partial f(x)$ if and only if $x\in\partial f^\bullet(w)$.
\medskip\\
{\rm (Proof) Note that both statements, $w\in \partial f(x)$ and 
$x\in\partial f^\bullet(w)$, are equivalent to that 
$f(x)+f^\bullet(w)\le \langle w,x\rangle$.
\pend}
\end{lemma}

\noindent
{\bf Remark}: 
When $f: \mathbb{R}^E\to\mathbb{R}\cup\{+\infty\}$ is an 
M${}^\natural$-convex function, subdifferentials $\partial f^\bullet(w)$ 
for all $w\in{\rm dom}(f^\bullet)$ are generalized polymatroids (or 
M${}^\natural$-convex sets). Furthermore, if $f$ is 
defined on integer lattice $\mathbb{Z}^E$, then $\partial f^\bullet(w)$
for all $w\in{\rm dom}(f^\bullet)$ are integral generalized polymatroids 
(restricted on $\mathbb{Z}^E$).
\pend\medskip

For more details about M-/M${}^\natural$-convex functions 
and L-/L${}^\natural$-convex functions see 
\cite{Murota03a,Murota2016} and \cite[Chapter VII]{Fuji05}.

\section{Compression of M${}^\natural$-convex Functions}\label{sec:compression}

In this section we introduce a new transformation, called {\it compression}, 
of M${}^\natural$-convex functions defined on the integer lattice 
$\mathbb{Z}^E$. 
The compression of an M${}^\natural$-convex function 
$f: \mathbb{Z}^E\to\mathbb{R}\cup\{+\infty\}$ is a transformation of 
the M${}^\natural$-convex function $f$ to an M-convex function 
$\hat{f}: \mathbb{Z}^E\to\mathbb{R}\cup\{+\infty\}$.

Consider any M${}^\natural$-convex function 
$f: \mathbb{Z}^E\to\mathbb{R}\cup\{+\infty\}$.
We suppose the following:
\begin{itemize}
\item The effective domain ${\rm dom}(f)$ is bounded and full-dimensional. 
That is, ${\rm dom}(f)$ is a full-dimensional generalized polymatroid 
${\rm P}(f^*,g^*)$ (with finite $f^*(E) > g^*(E)$).
\end{itemize}
For each integer $\alpha$ such that $f^*(E)\ge\alpha\ge g^*(E)$ let 
$f_{(\alpha)}: \mathbb{Z}^E\to\mathbb{R}\cup\{+\infty\}$ be the M-convex 
function defined by
\begin{equation}\label{eq:comp1}
 f_{(\alpha)}(x)=\left\{
 \begin{array}{ll}
   f(x) & {\rm if\ }\ x(E)=\alpha \\
   +\infty & {\rm otherwise}
 \end{array}
 \right. \quad (x\in\mathbb{Z}^E)
\end{equation}
(cf.~\cite{Murota03a,Murota2016,Murota2019}). 
We call $f_{(\alpha)}$ the {\it $\alpha$-section}
of $f$.
Then put 
\begin{equation}\label{eq:comp1a}
I_f=\{\alpha\in\mathbb{Z}\mid f^*(E)\ge\alpha\ge g^*(E)\}
\end{equation}
 and consider the {\sl convolution} of all the sections  $f_{(\alpha)}$ 
$(\alpha\in I_f)$,  which is given by
\begin{equation}\label{eq:comp2}
 \hat{f}(x)=\min\left\{\sum_{\alpha\in I_f}f_{(\alpha)}(y_\alpha)
              \ \Bigl|\  x=\sum_{\alpha\in I_f}y_\alpha, 
                 \ \forall \alpha\in I_f: y_\alpha\in\mathbb{Z}^E\right\}
 \quad (x\in\mathbb{Z}^E).
\end{equation}
where note that $f_{(\alpha)}(y_\alpha)<+\infty$ only if $y_\alpha(E)=\alpha$.
The convolution of M-convex functions $f_{(\alpha)}$ $(\alpha\in I_f)$ is 
an M-convex function (\cite{Murota03a,Murota2016,Murota2019}). 
We call $\hat{f}$ the {\it compression} of $f$. 

\begin{theorem}\label{th:comp1}
For the compression $\hat{f}$ of an M${}^\natural$-convex function 
$f: \mathbb{Z}^E\to\mathbb{R}\cup\{+\infty\}$ 
the Legendre-Fenchel transform of $\hat{f}$ is given by
\begin{equation}\label{eq:comp3}
 \hat{f}^\bullet = \sum_{\alpha\in I_f} f_{(\alpha)}{}^\bullet.
\end{equation}
We also have
\begin{equation}\label{eq:comp3a}
 {\rm dom}(\hat{f})=\sum_{\alpha\in I_f}{\rm dom}(f_{(\alpha)}),
\end{equation}
where the right-hand side is the Minkowski sum of the effective domains 
${\rm dom}(f_{(\alpha)})$ for all $\alpha\in I_f$.
\medskip\\
{\rm (Proof)
For any $w\in(\mathbb{R}^E)^*$ we have from (\ref{eq:comp2})
\begin{eqnarray}\label{eq:comp4}
 \hat{f}^\bullet(w)
  &=&\sup\{\langle w,x\rangle-\hat{f}(x)\mid x\in\mathbb{Z}^E\}\nonumber\\
&=&\sup\left\{\sum_{\alpha\in I_f}
       (\langle w,y_\alpha\rangle-f_{(\alpha)}(y_\alpha))
              \ \Bigl|\ \ \forall \alpha\in I_f: 
                y_\alpha\in\mathbb{Z}^E\right\}\nonumber\\
&=&\sum_{\alpha\in I_f}\sup\{\langle w,x\rangle-f_{(\alpha)}(x)
                      \mid x\in\mathbb{Z}^E\}\nonumber\\
&=&\sum_{\alpha\in I_f} f_{(\alpha)}{}^\bullet(w).
\end{eqnarray}
This implies (\ref{eq:comp3}) and (\ref{eq:comp3a}) 
because of  the definition of the Legendre-Fenchel transform. 
\pend}
\end{theorem}

\section{M${}^\natural$-convex Functions and Strip Decomposition}
\label{sec:strip}

Now, we introduce a concept of the {\it strip decomposition} of  
an M${}^\natural$-convex function defined on the integer lattice 
$\mathbb{Z}^E$. 
Let $f: \mathbb{Z}^E\to\mathbb{R}\cup\{+\infty\}$ be 
an M${}^\natural$-convex function with a full-dimensional and bounded 
effective domain ${\rm dom}(f)$. 

\subsection{Strip decomposition of M${}^\natural$-convex functions}
\label{sec:strip-decomposition}

For the compression $\hat{f}$ of $f$ and 
for any $w\in(\mathbb{R}^E)^*$ let us define 
$D(\hat{f},w)\subseteq \mathbb{Z}^E$ and 
$D(f_{(\alpha)},w)\subseteq \mathbb{Z}^E$ $(\alpha\in I_f)$ by 
\begin{equation}\label{eq:comp5}
 D(\hat{f},w)
 ={\rm Arg}\min\{\hat{f}(x)-\langle w,x\rangle\mid x\in\mathbb{Z}^E\} 
\end{equation}
\begin{equation}\label{eq:comp6}
 D(f_{(\alpha)},w)
={\rm Arg}\min\{f_{(\alpha)}(x)-\langle w,x\rangle\mid x\in\mathbb{Z}^E\} 
  \quad (\alpha\in I_f),
\end{equation}
where recall (\ref{eq:comp1a}) for the definition of $I_f$. 
We call $D(\hat{f},w)$ and $D(f_{(\alpha)},w)$ $(\alpha\in I_f)$ 
{\it linearity domains of} $\hat{f}$ and 
$f_{(\alpha)}$ $(\alpha\in I_f)$, respectively, {\it associated with} $w$.
We see from Lemma~\ref{lem:sbd1} that for every $x\in\mathbb{Z}^E$ 
we have 
\begin{equation}\label{eq:comp6a}
x\in D(\hat{f},w) \Longleftrightarrow w\in\partial \hat{f}(x) 
\Longleftrightarrow x\in \partial\hat{f}^\bullet(w).
\end{equation}
This means that the linearity domains of $\hat{f}$ are exactly 
the subdifferentials of $\hat{f}^\bullet$. 

Let $\mathbb{S}$ be the collection of all maximal linearity domains 
of $\hat{f}$ (or maximal subdifferentials of $\hat{f}^\bullet$). 
Then $\mathbb{S}$ gives a polyhedral division of ${\rm dom}(\hat{f})$ into 
generalized polymatroids. 
Since ${\rm dom}(f)$ is full-dimensional, the dimension of 
${\rm dom}(\hat{f})$ is equal to $|E|-1(=n-1)$.

\begin{figure}[h]
 \begin{center}
  \leavevmode
  \includegraphics[width=10cm,height=7cm,clip]{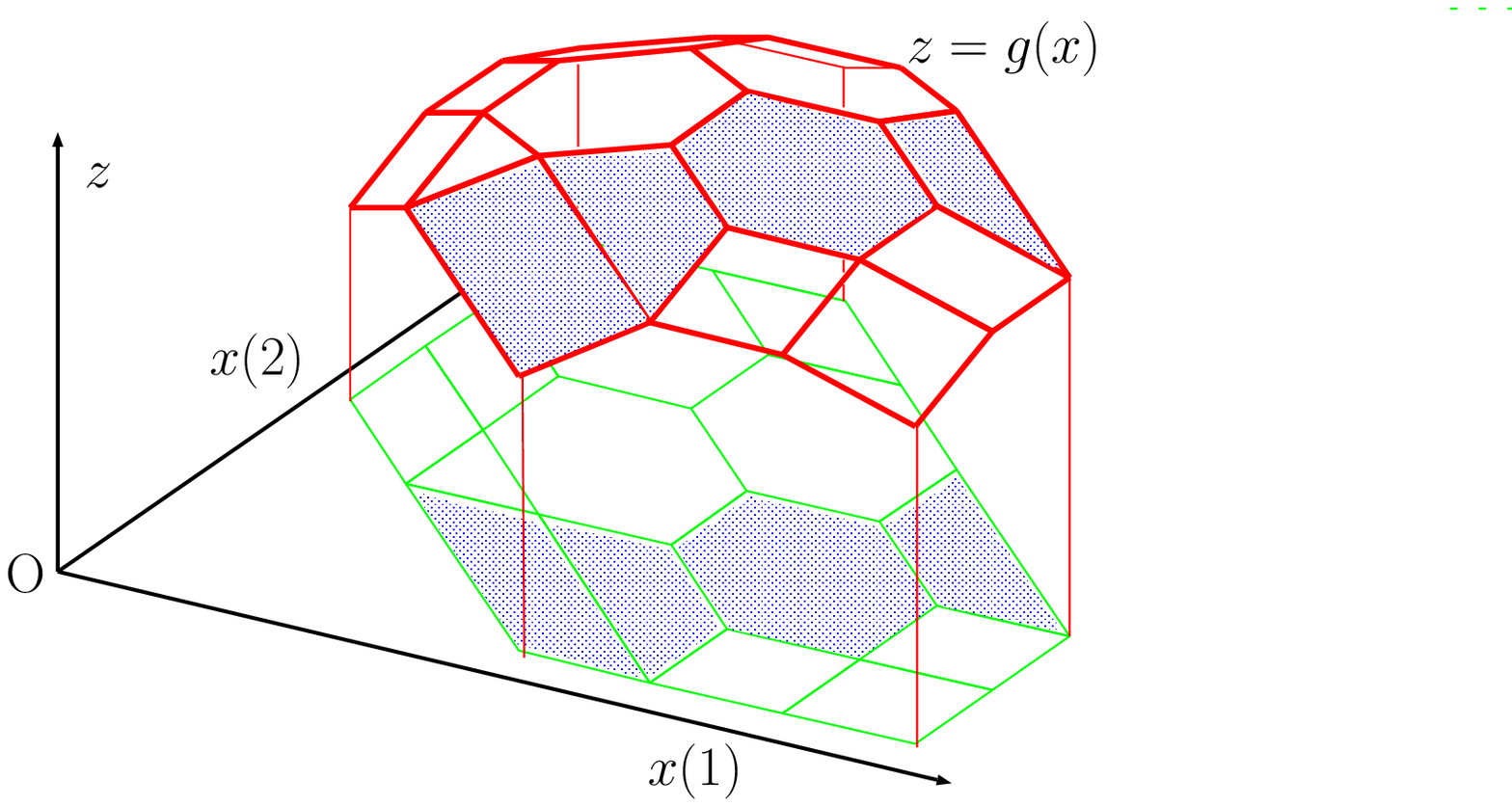}
 \end{center}
\caption{A strip of an M${}^\natural$-concave function $g$ 
indicated by shade.}
\label{fig:strip}
\end{figure}

For each maximal linearity domain $S\in\mathbb{S}$ of $\hat{f}$ 
let $w$ be a vector in $(\mathbb{R}^E)^*$ 
such that $S=\partial\hat{f}^\bullet(w)$ 
and put ${D}(\alpha,S)=D({f}_{(\alpha)},w)$ 
$(\alpha\in I_f)$, which 
are independent of the choice of $w$ satisfying 
$S=\partial\hat{f}^\bullet(w)$. 
Define ${D}(S)=\bigcup\{{D}(\alpha,S)\mid \alpha\in I_f\}$ and let
$f^{{D}(S)}$ be the restriction of $f$ on ${D}(S)$. 
Then we call $f^{{D}(S)}=(f_{(\alpha)}^{{D}(\alpha,S)}\mid \alpha\in I_f)$ 
a {\it strip of $f$ associated with $S\in\mathbb{S}$}, where 
$f_{(\alpha)}^{{D}(\alpha,S)}$ is
the restriction of $f_{(\alpha)}$ on ${D}(\alpha,S)$.
(See Figure~\ref{fig:strip} for an example of a strip of 
an M${}^\natural$-concave function.) 
The collection of the strips $f^{{D}(S)}$ for all $S\in\mathbb{S}$ 
is called the {\it strip decomposition} of $f$.

\subsection{Strips viewed from parametric optimization}\label{sec:parametric}

For any strip $f^{{D}(S)}=(f_{(\alpha)}^{{D}(\alpha,S)}
\mid \alpha\in I_f)$ 
 of an M${}^\natural$-convex function $f$ associated with $S\in\mathbb{S}$ 
let $w$ be a vector in $(\mathbb{R}^E)^*$ such that 
$S=\partial\hat{f}^\bullet(w)$. 
Then consider a parametric optimization problem
${\bf P}(\lambda)$ with a parameter $\lambda\in\mathbb{R}$ described 
as follows.
\begin{equation}\label{eq:par1}
 {\bf P}(\lambda):\ \ {\rm Minimize }\ f(x)-\langle w+\lambda{\bf 1},x\rangle\ 
    {\rm \ subject\ to\ }\ x\in{\rm dom}(f),
\end{equation}
where ${\bf 1}=\chi_E$ is the $n$-dimensional vector of all ones. 
We then have the following theorem.

\begin{theorem}\label{th:par1}
For $w\in(\mathbb{R}^E)^*$ chosen as above there exist a finite sequence
of values 
$\lambda_0=-\infty<\lambda_1<\cdots<\lambda_p<\lambda_{p+1}=+\infty$ and 
that of integers $k_0=g^*(E)<k_1<\cdots<k_p=f^*(E)$ such that 
the set $X^*(\lambda)$ of optimal solutions of ${\bf P}(\lambda)$ for each 
$\lambda\in\mathbb{R}$ is given by
\begin{equation}\label{eq:par2}
 X^*(\lambda)=\left\{
  \begin{array}{ll}
     \bigcup_{\alpha=k_{\ell-1}}^{\ k_\ell}{D}(\alpha,S)
         &  {\rm if\ } \lambda=\lambda_\ell \quad(\ell=1,\cdots,p)\\
     {D}(k_\ell,S) & 
         {\rm if\ }\ \lambda\in(\lambda_\ell,\lambda_{\ell+1}) 
                          \quad   (\ell=0,\cdots,p).
  \end{array}
  \right.
\end{equation}
{\rm (Proof)
Because of the discrete structure of the M${}^\natural$-convex 
function $f$ and the assumption that ${\rm dom}(f)$ is bounded, 
there exists a finite sequence of values 
$\lambda_0=-\infty<\lambda_1<\cdots<\lambda_p<\lambda_{p+1}=+\infty$ 
such that 
\begin{enumerate}
\item for each $i=0,1,\cdots,p$ Problems ${\bf P}(\lambda)$ for all
$\lambda\in(\lambda_i,\lambda_{i+1})$ have one and the same optimal 
solution set, and
\item the set $X^*(\lambda_i)$ of optimal solutions of ${\bf P}(\lambda_i)$
for each $i=1,\cdots,p$ consists of more than one optimal solution and 
we have  
$X^*(\lambda_i)\cap X^*(\lambda_{i+1})=X^*(\lambda)$ for each 
$i=1,\cdots,p-1$ and $\lambda\in(\lambda_i,\lambda_{i+1})$.
\end{enumerate}
Hence the optimal solution sets $X^*(\lambda)$ are expressed as (\ref{eq:par2})
for a sequence of some integers $k_0=g^*(E)<k_1<\cdots<k_p=f^*(E)$ 
that gives a division of the interval $I_f$.
\pend}
\end{theorem}
Here it should be noted that 
values $\lambda_i$ $(i=1,\cdots,p)$
depend on the choice of $w\in S$, while the 
vectors  
$w+\lambda_i{\bf 1}$ $(i=1,\cdots,p)$ are uniquely determined by $f$ and 
$S\in\mathbb{S}$, because of the assumptions that ${\rm dom}(f)$ is 
full-dimensional and $S\in\mathbb{S}$ is a maximal linearity domain of 
$\hat{f}$.

\section{Valuated Generalized Matroids}\label{sec:VGM}

In this section we further investigate the structures of strips
and their compressions for a class of valuated generalized matroids, 
which are special M${}^\natural$-convex functions defined on 
the unit hypercube $\{0,1\}^E$, in more details.

Let $f:\mathbb{Z}^E\to\mathbb{Z}\cup\{+\infty\}$ be an M${}^\natural$-convex 
function such that ${\rm dom}(f)=\{0,1\}^E$, which is called a {\it valuated
generalized matroid}. 
For any $X\subseteq E$ we often write $f(X)$ as $f(\chi_X)$ 
and regard $f$ as a function on $2^E$ in the sequel.

\subsection{Compression of valuated generalized matroids and 
valuated permutohedra}\label{sec:vflm}

For a valuated generalized matroid $f: 2^E\to\mathbb{Z}$ the compression 
$\hat{f}$ of $f$ given by (\ref{eq:comp2}) becomes
\begin{equation}\label{eq:vflm1}
 \hat{f}(x)=\min\left\{\sum_{\alpha\in [n]}f(Y_\alpha)
              \ \Bigl|\  x=\sum_{\alpha\in [n]}\chi_{Y_\alpha}, 
                 \ \forall \alpha\in [n]: Y_\alpha\in{E\choose\alpha}\right\}
 \quad (x\in\mathbb{Z}^E),
\end{equation}
where ${E\choose\alpha}=\{X\subseteq E\mid |X|=\alpha\}$ for $\alpha\in[n]$,
and we define $\hat{f}(x)=+\infty$ if the minimum on the right-hand side 
does not exist for $x\in\mathbb{Z}^E$. 
Then the effective domain of the compression $\hat{f}$ given by 
(\ref{eq:comp4}) is expressed by the following Minkowski sum:
\begin{equation}\label{eq:vflm2}
 {\rm dom}(\hat{f})
  =\sum_{\alpha\in[n]}\left\{\chi_Y\ \Bigl|\ Y\in{E\choose\alpha}\right\}.
\end{equation}
Recall that $E=[n]$.

\begin{theorem}\label{th:vflm1}
The 
effective domain ${\rm dom}(\hat{f})$ of the compression $\hat{f}$ 
is a permutohedron.
Hence the compression $\hat{f}$ is a valuated permutohedron, whose linearity
domains are sub-permutohedra.
\medskip\\
{\rm (Proof) The right-hand side of (\ref{eq:vflm2}) is the Minkowski 
sum of the sets of the characteristic vectors of bases of uniform matroids 
$U_{\alpha,n}$
of rank $\alpha$ for $\alpha\in[n]$. Hence it is a base polyhedron 
whose every extreme point 
(a greedy solution in the sense of Edmonds \cite{Edmonds70}) is a permutation 
$(\pi(1),\cdots,\pi(n))\in\mathbb{Z}^n$ of $[n]$ and {\sl vice versa}. 
It follows that (the convex hull of) 
${\rm dom}(\hat{f})$ is a permutohedron and $\hat{f}$ is 
a valuated permutohedron. 
Moreover, for any generic $w\in(\mathbb{R}^E)^*$ and every $\alpha\in[n]$, 
 $D(f_{(\alpha)},w)$ in
(\ref{eq:comp6}) is a singlton, $\chi_{F(\alpha)}$ say. Then, 
for each $\alpha=1,\cdots,n$  we have $F(\alpha-1)\subset F(\alpha)$ and 
$\chi_{F(\alpha)} - \chi_{F(\alpha-1)}=\chi_i$ for some $i\in[n]\ (=E)$ 
with $F(0)=\emptyset$. 
Hence sets $F(\alpha)$ for $\alpha=1,\cdots,n$\ \ form a complete flag
\begin{equation}\label{eq:vflm2a}
  \emptyset=F(0)\subset F(1)\subset F(2)\subset\cdots\subset F(n)=[n]
\end{equation}
and it determines a permutation $\pi$ of $[n]$ with the permutation vector 
$v_\pi=\sum_{\alpha=1}^n\chi_{F(\alpha)}$. 
It follows that every linearity domain 
of the compression $\hat{f}$ is a sub-permutohedron.
\pend}
\end{theorem}

\subsection{Strips of valuated generalized matroids and flag matroids}
\label{sec:strip-flag}

For every $S\in\mathbb{S}$ we have the strip 
$f^{{D}(S)}=(f_{(\alpha)}^{{D}(\alpha,S)}\mid \alpha=0,1,\cdots,n)$ 
of $f$ associated with $S\in\mathbb{S}$, which is characterized as follows.
We identify $\chi_X$ with $X$ for any $X\subseteq [n]$.

\begin{theorem}\label{th:vflm2}
For each $\alpha=0,1,\cdots,n$ we have a base family ${D}(\alpha,S)$ of 
a matroid $([n],\rho_\alpha^S)$ with a rank 
function $\rho_\alpha^S$ satisfying $\rho_\alpha^S([n])=\alpha$. 
Moreover, the sequence of $(\rho_\alpha^S\mid \alpha=0,1,\cdots,n)$ is that of 
strong maps, i.e., a flag matroid. 
\medskip\\
{\rm (Proof) The present theorem follows from the definition of the strip 
$f^{{D}(S)}=(f_{(\alpha)}^{{D}(\alpha,S)}\mid \alpha=0,1,\cdots,n)$ and 
the assumption that $f$ is a valuated generalized matroid.
\pend}
\end{theorem}

We call the flag matroid $(\rho_\alpha^S\mid \alpha=0,1,\cdots,n)$ the {\it 
flag matroid associated with a strip} $S\in\mathbb{S}$ 
(or a {\it flag-matroid strip}\/)
of the valuated generalized matroid $f$.
We see from Theorems~\ref{th:vflm1} and \ref{th:vflm2} the following.

\begin{theorem}\label{th:vfl32}
Every valuated generalized matroid $f$ induces a valuated permutohedron 
$\hat{f}$ 
by its compression and each flag-matroid strip of $f$ corresponds to a maximal 
linearity domain, a sub-permutohedron, of the induced valuated permutohedron.
\end{theorem}

Every valuated generalized matroid is regarded as a valuated permutohedron 
endowed with {\sl valuated} flag-matroid strips, one for each 
maximal linearity domain of it.

\section{Concluding Remarks}\label{sec:conclusion}

We have introduced the concepts of strip decomposition and compression 
of M${}^\natural$-convex functions. 
We have examined the structures of valuated generalized-matroids 
by considering the strip decomposition of a valuated generalized-matroid
into flag-matroid strips. The compression of a valuated generalized matroid 
induces a valuated permutohedron, a special M-convex function 
of Murota~\cite{Murota03a}. 
We thus have a new transformation, which we call the 
 compression, of a valuated generalized-matroid 
(an M${}^\natural$-convex function) 
to a valuated permutohedron (a special M-convex function).
Every Bruhat interval polytope is known to be a sub-permutohedron, 
due to Tsukerman and Williams~\cite{TsukermanWilliams2015}.
It is interesting to investigate Bruhat interval polytopes 
from a point of view of the strip decomposition of 
valuated generalized-matroids and also from a point of view 
of valuated permutohedra.

\section*{Acknowledgments}

Satoru Fujishige's work is supported by JSPS KAKENHI Grant Number~19K11839 
and Hiroshi Hirai's work is supported by JSPS KAKENHI Grant Number JP17K00029 
and JST PRESTO Grant Number JPMJPR192A, Japan.

\end{document}